\documentstyle[12pt]{article}
\def\P{{\bf P}}
\def\C{{\cal C}}
\def\K{{\cal K}}

\begin{document}

\begin{center}
{\Large \bf Base manifolds for fibrations of projective}\\
{\Large \bf  irreducible symplectic  manifolds}

\bigskip
 {\large \bf Jun-Muk Hwang}\footnote{This work was supported by the Korea Research Foundation Grant
funded by the Korean Government(MOEHRD)(KRF-2006-341-C00004).}

 \end{center}

\bigskip
\begin{abstract}
Given a projective irreducible symplectic manifold $M$ of dimension
$2n$, a projective manifold $X$  and  a surjective holomorphic map
$f:M \rightarrow X$ with connected fibers of positive dimension, we
prove that $X$ is biholomorphic to the projective space of dimension
$n$. The proof is obtained by exploiting two geometric structures at
general points of $X$: the affine structure arising from the action
variables of the Lagrangian fibration $f$ and the structure defined
by the variety of minimal rational tangents on the Fano manifold
$X$.
\end{abstract}


\medskip
2000MSC:14J40, 14J45

\bigskip
\section{Introduction}

We work in the category of complex analytic sets.  A connected
complex manifold $M$ of dimension $2n$ equipped with a holomorphic
symplectic from $\omega \in H^0(M, \Omega^2_M)$ is called a {\it
holomorphic symplectic manifold}. A subvariety $V$ of $M$ is said to
be {\it Lagrangian} if $V$ has dimension $n$ and the restriction of
$\omega$ on the smooth part of $V$ is identically zero. A simply
connected projective algebraic manifold $M$ is called a {\it
projective irreducible symplectic manifold} if $M$ has a symplectic
form $\omega$ such that $H^0(M, \Omega^2_M) = {\bf C} \omega.$ It is
remarkable that fibrations of projective irreducible
 symplectic manifolds are of very special form, as described in the
following theorem  due to D. Matsushita.

\medskip
{\bf Theorem 1.1} {\it Let $M$ be a projective irreducible
symplectic manifold of dimension $2n$. For a projective manifold $X$
and  a surjective holomorphic map $f:M \rightarrow X$ with connected
fibers of positive dimension, the following holds.

(1) $X$ is a Fano manifold of dimension $n$ with Picard number 1.

(2) A general fiber of $f$ is biholomorphic to an abelian variety.

(3) The underlying subvariety of every fiber of $f$ is Lagrangian.

(4) All even Betti numbers of $X$ are equal to 1 and all odd Betti
numbers of $X$ are equal to 0.}

\medskip
(1), (2) and (3)  in Theorem 1.1 were proved in [Ma1] and [Ma2].
These results led to the question whether the base manifold $X$ is
the complex projective space (cf. [Hu, 21.4]). The result of [Ma3]
verifies Theorem 1.1 (4), i.e.,  that the Betti numbers of $X$ are
indeed equal to those of $\P_n$. Some special cases of this question
were studied in [CMS, Section 7] and [Ng].

\medskip
Our goal is to give an affirmative answer to the question as
follows.

\medskip
{\bf Theorem 1.2} {\it In the setting of Theorem 1.1,  $X$ is
biholomorphic to $\P_n$.}

\medskip
There are two geometric ingredients in the proof of Theorem 1.2: the
theory of varieties of minimal rational tangents and the theory of
Lagrangian fibrations. On the one hand, the theory of varieties of
minimal rational tangents describes a certain geometric structure
arising from minimal rational curves at general points of a Fano
manifold $X$ with $b_2(X)=1$ (cf. [HwMo1], [HwMo2]). This geometric
structure has differential geometric properties reflecting special
features of the deformation theory of minimal rational curves. On
the other hand, the theory of Lagrangian fibrations, or
equivalently, the theory of completely integrable Hamiltonian
systems, provides an affine structure at general points of the base
manifold $X$ via the classical action variables (cf. [GuSt, Section
44]). Our strategy to prove Theorem 1.2 is to exploit the interplay
of these two geometric structures on the base manifold $X$. Under
the assumption that $X$ is different from $\P_n$, the condition
$b_2(X)=1$ forces the geometric structure arising from the variety
of minimal rational tangents to be `non-flat', while the affine
structure arising from the action variables is naturally `flat'.
These two structures interact via the monodromy of the Lagrangian
fibration, leading to a contradiction. To be precise, two separate
arguments are needed depending on whether the dimension $p$ of the
variety of minimal rational tangents is positive or zero. The easier
case of $p>0$ is handled by a topological argument using $b_4(X)
=1$. The more difficult case of $p=0$ needs a deeper argument,
depending on the local differential geometry of the variety of
minimal rational tangents.

\medskip
It is expected that an analog of Theorem 1.2 holds under the weaker
assumption that $M$ is a K\"ahler irreducible symplectic manifold
and $X$ is a compact K\"ahler manifold (cf. [Hu,24.8]). As far as we
see, the only point where our proof of Theorem 1.2 uses the
projectivity assumption in an essential way is the result $b_4(X)
=1$ from [Ma3], which is used to handle the case of $p>0$. Thus if
[Ma3] is generalizable to K\"ahler setting, so is our result. On the
other hand, after seeing the first draft of this paper, D.
Matsushita informed us that he has a new approach to the case of
$p>0$, which may lead to the generalization to K\"ahler setting,
even allowing some singularity on $M$.

\medskip
It is also expected that an analog of Theorem 1.2 holds under the
weaker assumption that $X$ is a normal projective variety. In fact,
under this weaker assumption, [Ma1] showed that the singularity of
$X$ is mild and $X$ must be a Fano variety of Picard number 1.
However, our approach to Theorem 1.2 uses the smoothness of $X$ in a
crucial way and seems difficult to generalize to the singular
setting.

\medskip
In Section 2, we present  results about varieties of minimal
rational tangents of Fano manifolds that are independent of
Lagrangian fibrations. In Section 3, we present results about
Lagrangian fibrations that are independent of Fano manifolds.
These two theories are played against each other to prove Theorem
1.2, first in the case of $p>0$ in Section 4 and then in the case
of $p=0$ in Section 5.

\section{Results on varieties of minimal rational
tangents}

In this section, we present several results about varieties of
minimal rational tangents, especially when they are linear. Most
of these are already known, for which we only give brief
explanation with precise references.

\medskip Throughout this section, we will denote by $X$
 an $n$-dimensional  Fano manifold with $b_2(X)=1$. An irreducible
component ${\cal K}$ of the space of rational curves on $X$ is
called a {\it minimal component} if for a general point $x \in X$,
the subscheme ${\cal K}_x$ of ${\cal K}$ consisting of members
passing through $x$ is non-empty and complete. In this case, the
subvariety ${\cal C}_x$ of the projectivized tangent space
${\P}T_x(X)$ consisting of the tangent directions at $x$ of members
of ${\cal K}_x$ is called the {\it variety of minimal rational
tangents} at $x$ (see [HwMo2] for more details). For a general
member $C$ of ${\cal K}$, the normalization $\nu: \hat{C}
\rightarrow C \subset X$ is an immersion of $\P_1$.

\medskip
 The following
result is proved in [HwRa, Corollary 2.2].

\medskip
{\bf Proposition 2.1} {\it Given a minimal component ${\cal K}$ on
$X$ and a general point $x \in X$,  let ${\cal C}_1$ be a
component of ${\cal C}_x$. Denote by ${\cal K}_1$ the
corresponding component of ${\cal K}_x$.  For a general member $C$
of ${\cal K}_1$, let $S_C \subset \P T^*_x(X)$ be the linear
subspace of the projectivized cotangent space that is the image of
the evaluation of $H^0(\hat{C}, \nu^*T^*(X))$ at $x$. Then the
closure of the union of $S_C$ as $C$ varies over general points of
${\cal K}_1$ is the dual variety of ${\cal C}_1 \subset \P
T_x(X)$.}

\medskip
Proposition 2.1 is useful when combined with the following.

\medskip
{\bf Proposition 2.2} {\it Given a minimal component ${\cal K}$ on
$X$ and a general point $x \in X$, suppose the dual variety ${\cal
C}_1^* \subset {\bf P}T^*_x(X)$ of a component ${\cal C}_1$ of
${\cal C}_x$ is linearly degenerate, i.e., contained in a
hyperplane in $\P T^*_x(X)$. Then the component ${\cal C}_1$ is a
linear subspace of ${\bf P}T_x(X)$.}

\medskip
{\it Proof}.  If ${\cal C}_1^*$ is linearly degenerate in ${\bf
P}T^*_x(X)$, then ${\cal C}_1$ is a cone. Thus Proposition 2.2 is
equivalent to  [HwMo2, Proposition 13], which says that ${\cal C}_1$
cannot be a cone unless it is a linear subspace. $\Box$

\medskip
The next proposition is [HwMo1, Lemma 4.2].

\medskip
{\bf Proposition 2.3} {\it Let $X'$ be a projective algebraic
manifold and $X$ be a Fano manifold with $b_2(X)=1$. Fix a minimal
component ${\cal K}$ on $X$.  Given a generically finite morphism
$\mu: X' \rightarrow X$ which is not birational, and given a
general member $C \subset X$ of ${\cal K}$, there exists a
component $C'$ of $\mu^{-1}(C)$ such that the restriction
$\mu|_{C'} : C' \rightarrow C$ is finite of degree $>1$.}

\medskip
For the rest of this section, we will make the following
assumption.

\medskip
({\it Assumption})  $X$ is an $n$-dimensional Fano manifold with
$b_2=1$, different from $\P_n$, and for some choice of ${\cal K}$,
the variety of minimal rational tangents at a general point  is
the union of linear subspaces of dimension $p \geq 0$.

\medskip
The condition $ X \neq \P_n$ implies that $p<n-1$ by [CMS, Theorem
0.2] and the condition $b_2(X)=1$ implies that ${\cal C}_x$ has at
least two irreducible components.
 In an analytic local neighborhood of $x$, each component of $\C_x$
defines a distribution. It is easy to see that this distribution
is integrable and the leaf through $x$ is an immersed $\P_{p+1}$.
More precisely, we have the following result from [Hw, Proposition
1].

\medskip
{\bf Proposition 2.4} {\it When $X$ satisfies the above
assumption, there exists a projective algebraic manifold $X'$ with
a generically finite holomorphic map $\mu:X' \rightarrow X$ of
degree $>1$ and a proper holomorphic map $\rho: X' \rightarrow Z$
onto a positive-dimensional projective manifold $Z$  such that
$\rho$ is a $\P_{p+1}$-bundle over a Zariski open subset
$Z_o\subset Z$, $\mu$ is unramified on $\rho^{-1}(Z_o)$, and each
member of $\K_x$ for a general $x \in X$ is the image of a line in
some ${\bf P}_{p+1}$-fibers of $\rho$. Here, $p + \dim Z = n-1$
and $\dim Z = \dim H^0(\hat{C}, \nu^*T^*(X))$ for the
normalization $\nu: \hat{C} \rightarrow C \subset X$ of a general
member $C$ of ${\cal K}$.  Moreover, if we set
 $P_{\zeta}:= \mu(\rho^{-1}(\zeta))$ for each $\zeta \in Z_o$,
 then the following holds.

  (a)
$P_{\zeta}$ is an immersed submanifold with trivial normal bundle
in $X$.

(b) $\mu|_{\rho^{-1}(\zeta)}$ is the normalization of $P_{\zeta}$.

 (c) For two distinct points $\zeta_1 \neq \zeta_2 \in Z_o$, the two
subvarieties $P_{\zeta_1}$ and $P_{\zeta_2}$ are distinct.}

\medskip
Here, when $X$ is a complex manifold and $V \subset X$ is a
subvariety, we  say that $V$ is an immersed submanifold if the
normalization $\hat{V}$ is smooth and the normalization map
$\nu:\hat{V} \rightarrow V \subset X$ is an immersion. The normal
bundle of $V$ means the vector bundle  on $\hat{V}$ defined as the
quotient of $\nu^*T(X)$ by the image of $T(\hat{V})$. If the
normal bundle of $V$ is a trivial bundle, we say that $V$ is an
immersed submanifold with trivial normal bundle.

\medskip
The first sentence of the following proposition is exactly [Hw,
Proposition 2.2], from which the second sentence follows because
$P_{\zeta}, \zeta \in Z_o$, is immersed with  trivial normal
bundle.

\medskip
{\bf Proposition 2.5} {\it In the setting of Proposition 2.4, for
a general point $x \in X$, the variety of minimal rational
tangents ${\cal C}_x$ is a disjoint union of linear subspaces in
${\bf P}T_x(X)$, i.e., for each $\zeta_1 \neq \zeta_2 \in Z_o$
with $x \in P_{\zeta_1} \cap P_{\zeta_2}$, the intersection of the
tangent spaces $T_x(P_{\zeta_1})$ and $T_x(P_{\zeta_2})$ in
$T_x(X)$ is zero. Consequently, the elements of $T_x^*(X)$
obtained by  evaluating $H^0(\rho^{-1}(\zeta_1), \mu^*T^*(X))$ and
$H^0(\rho^{-1}(\zeta_2), \mu^*T^*(X))$ at $x$ span $T^*_x(X)$. }

\medskip
In the situation of Proposition 2.4, it is convenient to introduce
the following notion, which already played an essential role in
[Hw]. We say that ${\cal K}$ is {\it multivalent} on an
irreducible hypersurface $H \subset X$, if the following holds:
when $H' \subset X'$ denotes the union of the components of
$\mu^{-1}(H)$ that are dominant over $H$ by $\mu$ and dominant
over $Z$ by $\rho$, the dominant map $H' \rightarrow H$ has degree
$>1$. The next proposition is a variation of [Hw, Proposition
3.2].

\medskip
{\bf Proposition 2.6} {\it In the setting of Proposition 2.4,
suppose that ${\cal K}$ is multivalent on an irreducible
hypersurface $H \subset X$. Then a general point $y \in H$ has two
open neighborhoods $W \subset W_0 $  satisfying the following.

(1) The fundamental group $\pi_1(W_0 \setminus H)$ is cyclic.

(2) For each general point $x \in W \setminus H$,  there are  two
distinct points $\zeta_1, \zeta_2 \in Z_o$ with $x \in P_{\zeta_1}
\cap P_{\zeta_2}$ and  $T_x(P_{\zeta_1}) \cap T_x(P_{\zeta_2}) = 0$
such that there exists a loop $\gamma_1$ (resp. $\gamma_2$) based at
$x$ lying on the smooth locus of $P_{\zeta_1} \setminus H$ (resp.
$P_{\zeta_2} \setminus H$), the homotopy class of which generates
the cyclic group $\pi_1(W_0 \setminus H, x)$.}

\medskip
{\it Proof}. The multivalence  assumption on $H$ implies that there
are two distinct  $P_{\zeta_3}$ and $P_{\zeta_4}$ passing through
$y$. The desired $P_{\zeta_1}$ and $P_{\zeta_2}$ can be  obtained as
small deformations of $P_{\zeta_3}$ and $P_{\zeta_4}$. Let us make
this more precise.

 Since ${\cal K}$ is multivalent on $H$, there
exist two distinct points $y_1, y_2 $ in $\mu^{-1}(y) \cap H'$ such
that the map $\mu$ is unramified at $y_1$ and $y_2$. We can choose
open neighborhoods $W_1 \subset X'$ of $y_1$, $W_2 \subset X'$ of
$y_2$ and $W_0 \subset X$ of $y$ with the following properties:

(a) $W_0$ is biholomorphic to the polydisc $\Delta^n$ and $W_0
\cap H$ corresponds to the coordinate hyperplane $\Delta^{n-1}
\subset \Delta^n.$

(b) $\mu(W_1) = \mu(W_2) = W_0$

 (c) $\mu|_{W_1}$ and $\mu|_{W_2}$ are biholomorphic over $W_0$.

 (d) $W_1 \cap H'$ and $W_2 \cap H'$ are
transversal to fibers of $\rho$ at the intersection points, i.e.,
their scheme-theoretic intersections with the fibers of $\rho$ are
smooth.

We can assume that $\rho$ is a ${\bf P}_{p+1}$-bundle near
$\rho(y_1)$ (resp. $\rho(y_2)$), so there exists an open
neighborhood $W'_1 \subset W_1$ of $y_1$ (resp. $W'_2 \subset W_2$
of $y_2$) such that for any $w \in W'_1$ (resp. $w \in W'_2$), $
\rho^{-1}(\rho(w)) \cap W'_1 $ (resp. $ \rho^{-1}(\rho(w)) \cap
W'_2$ ) is irreducible. Now choose $W$ to be a  neighborhood of
$y$ in $ \mu(W'_1) \cap \mu(W'_2)$.  Then (1) is automatic and (2)
can be seen from the following lemma. $\Box$

\medskip
{\bf Lemma 2.7} {\it Let $B = \Delta^n$ and $H \subset B$ be the
coordinate hyperplane $\Delta^{n-1} \subset \Delta^n$. Let $y \in
H$ be a point and $V \subset B$ be an irreducible  closed immersed
submanifold in $B$ such that $y \in V$ and an irreducible
component of the germ of $V$ at $y$ is non-singular and intersects
$H$ transversally at $y$. Then for any non-singular point $x$ of
$V$, there exists a loop based at $x$ lying on the smooth locus of
 $V \setminus H$ which generates the cyclic fundamental group
$\pi_1(B \setminus H, x)$.}

\medskip
{\it Proof}. Since $V$ is irreducible, it suffices to show this
for some non-singular point $x \in V \setminus H$. But this is
obvious if $x$ lies on the irreducible component of the germ of
$V$ intersecting $H$ transversally at $x$. $\Box$

\medskip
The next proposition is precisely [Hw, Proposition 3.1].

\medskip
{\bf Proposition 2.8} {\it Let $X$ be a Fano manifold with $b_2(X)
= b_4(X) =1$ satisfying the assumption of Proposition 2.4 with
$p>0$.  Then the minimal component ${\cal K}$ is multivalent on
every irreducible hypersurface $H \subset X$. }

\medskip
The condition $p>0$ in  Proposition 2.8 is crucial. In fact, when
$X$ is a Fano 3-fold defined by a linear section of the
6-dimensional Grassmannian of rank 3, there is a surface $H
\subset X$ such that the minimal component ${\cal K}$ is not
multivalent on $H$.

\section{ Results on Lagrangian fibrations}

In this section, we present some results about Lagrangian
fibrations.

\medskip
In this paper, the phrase   `Lagrangian fibration' will have the
following restrictive meaning. Let $(M, \omega)$ be a holomorphic
symplectic manifold  and $B$ be a complex manifold of dimension $n$.
A proper surjective holomorphic map $f:M \rightarrow B$ is  a {\it
Lagrangian fibration} if it satisfies the statements (2) and (3) of
Theorem 1.1, i.e., if a general fiber of $f$ is biholomorphic to an
abelian variety and the underlying subvariety of every fiber of $f$
is Lagrangian.

Given a Lagrangian fibration $f:M \rightarrow B$, the locus $D
\subset B$ of the critical values of $f$ is a hypersurface (e.g.
[HO, Proposition 3.1]), if non-empty. We will call $D$ the {\it
critical set} of $f$. The following is well-known, see, for
example, [HO, Proposition 3.2] for a proof.

\medskip
{\bf Proposition 3.1} {\it Let $f: M \rightarrow B$ be a
Lagrangian fibration. Given $b \in B$, if the underlying reduced
variety of $f^{-1}(b)$ is an abelian variety, then $f^{-1}(b)$ is
smooth, i.e., $b$ is not in the critical set of $f$.}

\medskip
A proper holomorphic map  $f: M \rightarrow B$ between two
connected complex manifolds is called a {\it smooth abelian
fibration} if every fiber is biholomorphic to an abelian variety.
Recall that for a smooth abelian fibration, a choice of a base
point $b \in B$ gives rise to the monodromy representation
$$\pi_1(B, b) \longrightarrow {\rm GL}(H_1(f^{-1}(b), {\bf Z})).$$

The following result is due to Matsushita [Ma4]. Since its proof in
[Ma4] is buried in a longer argument to prove a much more
substantial result, we will give a sketch of the proof for readers'
convenience.

\medskip
{\bf Proposition 3.2} {\it Let $f: M \rightarrow B$ be a
Lagrangian fibration with a non-empty critical set $D \subset B$.
  Assume that $B$ is biholomorphic to the polydisc
$\Delta^n$ and the critical set $D \subset B$ as a subvariety is
biholomorphic to a coordinate hyperplane $\Delta^{n-1} \subset
\Delta^n$. Assume furthermore that there exists a line bundle on
$M$ which is ample on each fiber. Then the monodromy of the smooth
abelian fibration $ M\setminus f^{-1}(D) \rightarrow B \setminus
D$ is non-trivial, i.e., for any point $b \in B \setminus D$, the
image of the representation
$$\pi_1( B \setminus D, b) \longrightarrow {\rm GL}(H_1(f^{-1}(b),
{\bf Z}))$$  is not the identity.}

\medskip
{\it Proof}. Assume that the monodromy is trivial. We can choose a
cyclic cover of $B$ along $D$, $\xi: B' \rightarrow B$ such that
the induced  fibration $\xi^*f: \xi^*M \rightarrow B'$ has a
section. If $f$ already has a section, we choose $\xi$ to be the
 identity map of $B$. Let $D' = \xi^{-1}(D)_{\rm red}$ be the ramification
 set of $\xi$. Using the triviality of the monodromy of $\xi^*f$ and
the ample line bundle, we can define the period map from $B'
\setminus D'$ to the Siegel upper half space, which extends to a
holomorphic map from $B'$ to the Siegel upper half space. Using
the extended period map, we can construct a smooth abelian
fibration $\tilde{f}: \tilde{M} \rightarrow B'$ with a
biholomorphic map
$$\Phi: \xi^*M \setminus [(\mu^*f)^{-1}(D')] \longrightarrow
\tilde{M} \setminus \tilde{f}^{-1}(D').$$ The existence of the
section for $\xi^*f$ ensures that $\Phi$ defines a bimeromorphic map
between $\xi^* M$ and $\tilde{M}$, by the argument of Nakayama [Nk,
Proposition 1.6].  The action of the Galois group $\Gamma$ of the
cyclic cover $\xi$ on $\xi^*M$ induces an action of $\Gamma$ on
$\tilde{M}$. If $\Gamma$ is not trivial, i.e., if $f$ doesn't have a
section, then  $\Gamma$ acts non-trivially on the fibers of
$\tilde{f}$ over $D'$. On the other hand, $\Gamma$ acts trivially on
the singular homology group of the fibers of $\tilde{f}$ over $D'$
by the triviality of the monodromy. This means that $\Gamma$ acts as
a translation on the fibers of $\tilde{f}$ over $D'$. Thus the
quotient of $\tilde{M}$ by $\Gamma$ is an abelian fibration
$\tilde{M}/\Gamma \rightarrow B$, the underlying reduced variety of
each fiber of which is an abelian variety. In particular,
$\tilde{M}/\Gamma$ contains no rational curve. The bimeromorphic map
$\Phi$ descends to a bimeromorphic map between $M$ and
$\tilde{M}/\Gamma$. Since $\hat{M}/\Gamma$ contains no rational
curve, $M$ and $\hat{M}$ are biholomorphic. By Proposition 3.1, the
critical set of $f$ is empty, a contradiction. $\Box$

\medskip
 Given a holomorphic symplectic manifold $(M, \omega)$, the
contraction with $\omega$ induces a natural isomorphism
$$\iota_{\omega}: T^*(M) \longrightarrow T(M).$$  For a Lagrangian fibration
$f:M \rightarrow B$, any $b \in B$ and any $z\in f^{-1}(b)$, we
have a homomorphism $$\iota_{\omega} \circ f^*: T_b^*(B)
\longrightarrow T_z(M).$$ If $z$ is a non-singular point of
$f^{-1}(b)$,  this induces an isomorphism
$$\iota_{b,z}: T^*_b(B) \longrightarrow T_z(f^{-1}(b)).$$   If $b$ is not in
the critical set, denote by ${\rm Aut}_o(f^{-1}(b))$ the identity
component of the automorphism group of the abelian variety
$f^{-1}(b).$ Then for any $z \in f^{-1}(b)$, $\iota_{b,z}$ induces
a natural unramified surjective group homomorphism
$$h_b: T^*_b(B) \longrightarrow {\rm Aut}_o(f^{-1}(b)).$$
An analog of this exists for arbitrary $b \in B$ as follows.

\medskip
{\bf Proposition 3.3} {\it Let $f: M \to B$ be a Lagrangian
fibration. For each point $b \in B$, there exists a canonical
homomorphism of complex Lie groups $$h_b: T^*_b(B) \longrightarrow
{\rm Aut}_o(f^{-1}(b)_{\rm red}),$$ where $T^*_b(B)$ is the vector
group and the target is the identity component of the automorphism
group of the underlying reduced variety of the fiber at $b$. }

\medskip
{\it Proof}. For each $v \in T^*_b(B)$, let $\tilde{v}$ be an
exact 1-form in a local analytic neighborhood $U$ of $b$ which
coincides with $v$ at $b$. The vector field $\iota_{\omega}(f^*
\tilde{v})$ on $f^{-1}(U)$ is a Hamiltonian vector field and is
tangent to the fibers of $f$. Thus it induces a derivation
$\vec{v}$ on $f^{-1}(b)_{\rm red}$. At a  point $z$ of
$f^{-1}(b)_{\rm red}$, $\vec{v}$ is just $\iota_{\omega} \circ
f^*(v)$. Thus this derivation is independent of the choice of
$\tilde{v}$ and depends only on $v$. By integrating the derivation
(e.g. [Ka, p.83, Korollar]), we get a natural group homomorphism
${\mathbf C} v \to {\rm Aut}_o (f^{-1}(b)_{\rm red})$. This
defines $h_b: T^*_b(B) \to {\rm Aut}_o (f^{-1}(b)_{\rm red})$ in a
canonical way. $\Box$

\medskip
{\bf Proposition 3.4} {\it Let $f: M \rightarrow B$ be a
Lagrangian fibration and let $V$ be an irreducible complete
variety in $B$. Then the following holds.

(i) There exists a canonical homomorphism of complex Lie groups
$$h_V: H^0(V, T^*(B)) \longrightarrow {\rm Aut}^f_o(f^{-1}(V)_{\rm red})$$ where
the target denotes the identity component of the complex Lie group
of automorphisms preserving the $f$-fibers of the underlying
reduced variety  of $f^{-1}(V)$.

 (ii) If $V$ is not contained in
the critical set of $f$, the group ${\rm Aut}^f_o(f^{-1}(V)_{\rm
red})$ is an abelian variety, which we denote by $A_V$.

(iii) If $V$ is not contained in the critical set of $f$, but has
a non-empty intersection with  the critical set, then the abelian
variety $A_V$  has dimension $\leq n-1$.

(iv) Suppose that $b \in V$ is not in the critical set and $z \in
f^{-1}(b)$. Let $T_e(A_V)$ be the tangent space of the abelian
variety $A_V$ in (ii)  at the identity. The homomorphism $h_V$ in
(i) induces
$$dh_V: H^0(V, T^*(B)) \longrightarrow T_e(A_V)$$ and there is a
commutative diagram
$$ \begin{array}{ccc} H^0(V, T^*(B)) & \longrightarrow & T_b^*(B) \\ d h_V \downarrow
& & \iota_{b, z} \downarrow \\ T_e(A_V) & \longrightarrow &
T_z(f^{-1}(b))
\end{array}$$ where the upper horizontal arrow is the evaluation
at $b$ and the lower horizontal arrow is the tangent map for the
$A_V$-orbit of $z$. }

\medskip
{\it Proof}. For each $\xi \in H^0(V, T^*(B))$, the automorphism
$h_V(\xi)$ is defined such that it acts on $y \in f^{-1}(V)_{\rm
red}$  by
$$ h_V(\xi) \cdot y = h_b( \xi_b) \cdot y$$ where $b = f(y) \in V$,  $h_b$ is the
homomorphism defined in Proposition 3.3 and $\xi_b$ is the value
of $\xi$ at $b$. This gives (i).  To see (ii), suppose that the
algebraic group ${\rm Aut}^f_o(f^{-1}(V)_{\rm red})$ is not an
abelian variety. Then it contains a linear algebraic subgroup by
Chevalley decomposition. This linear algebraic group must act
non-trivially on $f^{-1}(b)$ for a general point $b \in V$. But
linear algebraic groups cannot act on abelian varieties
nontrivially, a contradiction. Regarding (iii), if the dimension
of ${\rm Aut}^f_o(f^{-1}(V)_{\rm red})$ is $\geq n$, the fiber of
$f$ over each point of $V$ contains an orbit, which must be an
abelian variety of dimension $\geq n$. Thus each fiber is smooth
by Proposition 3.1, a contradiction to the assumption in (iii).
Finally, (iv) is immediate from the definition. $\Box$

\medskip
Recall that the cotangent bundle of a complex manifold has a
canonical holomorphic symplectic form on it. The following is a
geometric version of the classical action-angle variables for an
integrable Hamiltonian system. The proof in [GS, Theorem 44.2]
works verbatim in the holomorphic setting.

\medskip
{\bf Proposition 3.5} {\it Let $f:M \rightarrow B$ be a smooth
Lagrangian fibration with a Lagrangian section $\Sigma \subset M$.
The orbit of $\Sigma$ under the action of $T^*_b(B), b \in B,$ in
Proposition 3.3 gives rise to an unramified surjective holomorphic
map $\chi: T^*(B) \rightarrow M$ which commutes with the
projection to $B$ and the map $f$. Then the pull-back of the
symplectic form $\omega$ on $M$ by  $\chi$ coincides with the
standard symplectic structure of $T^*(B)$.}

\medskip
An immediate consequence is the following.

\medskip
{\bf Proposition 3.6} {\it Let $f:M \rightarrow B$ be a smooth
Lagrangian fibration with a Lagrangian section $\Sigma \subset M$.
Let $M' \subset M$ be a family of abelian subvarieties  of $f$,
i.e., a submanifold containing $\Sigma$ such that $f|_{M'}$ is a
smooth abelian fibration. Let ${\cal V} \subset T^*(B)$ be the
subbundle defined by $\chi^{-1}(M')$. Then ${\cal V}$, regarded as a
differential system on $B$, is involutive, i.e., the distribution
${\cal V}^{\perp} \subset T(B)$ annihilated by ${\cal V}$ is
integrable.}

\medskip
{\it Proof}. By Proposition 3.5, the submanifold $\chi^{-1}(\Sigma)
\subset T^*(B)$, i.e., the family of lattices for the family of
abelian varieties, is Lagrangian in $T^*(B)$. The family of lattices
for the family of abelian subvarieties  $M'$ is ${\cal V} \cap
\chi^{-1}(\Sigma)$. Thus the subbundle ${\cal V}$ is locally
generated by vectors lying in $\chi^{-1}(\Sigma)$. But Lagrangian
sections of $T^*(B)$ are just closed 1-forms on $B$. Thus the system
${\cal V}$ is locally generated by closed 1-forms, and consequently
it is involutive. $\Box$

\medskip
To use Proposition 3.6, we need to have a family of  abelian
subvarieties of a Lagrangian fibration. One construction leading to
a family of  abelian subvarieties is the following.

\medskip
{\bf Proposition 3.7} {\it Let $f: M \rightarrow B$ be a Lagrangian
fibration with non-empty critical set $D \subset B$. Suppose that
there exists a smooth proper morphism  $\rho: B \rightarrow Z$ of
relative dimension 1 onto a complex manifold $Z$. Assume that each
fiber of $\rho$ intersects the critical set $D$. Then for each fiber
$C$ of $\rho$, the abelian variety $A_C$ in the notation of
Proposition 3.4 is of dimension $n-1$. In particular, if we choose
an open subset $W \subset B\setminus D$ such that the restriction of
$f$ to $f^{-1}(W)$, $f^{-1}(W) \rightarrow W,$ admits a Lagrangian
section, then we have a family of $(n-1)$-dimensional abelian
subvarieties in the smooth Lagrangian fibration $f^{-1}(W)
\rightarrow W$, defined as the orbit of the Lagrangian section. The
distribution on $W$ induced by this family of abelian subvarieties
in the sense of Proposition 3.6 is just the fibers of $\rho $ on
$W$.  }

\medskip
{\it Proof}. For each $C$, $H^0(C, T^*(B))$ contains the
$(n-1)$-dimensional subspace corresponding to $\rho^*
T^*_{\rho(C)}(Z)$. This subspace injects into $T_e(A_C)$ by the
construction of the homomorphism $h_C$ in Proposition 3.4. Thus
$A_C$ is an abelian variety of dimension $\geq n-1$. By
Proposition 3.4 (iii), $A_C$ is an $(n-1)$-dimensional abelian
variety. In particular, the homomorphism $h_C$ sends
$T^*_{\rho(C)}(Z)$ onto $A_C$, which proves the last sentence of
Proposition 3.7. $\Box$

\section{Varieties of minimal rational
tangents for the base manifold}

\medskip
Now we go to the situation of Theorem 1.2. To start with, we recall
the following result, which can be proved in various ways. The proof
here was suggested by K. Oguiso.

\medskip
{\bf Proposition 4.1} {\it Let $f: M \rightarrow X$ be as in
Theorem 1.2. Then the critical set  $D$ is a non-empty
hypersurface.}

\medskip
{\it Proof}. We already mentioned in Section 3 that $D$ is a
hypersurface if it is non-empty. Suppose $D$ is empty. Then $f:M
\rightarrow X$ is a proper smooth morphism. Denote by ${\bf Z}_M$
and ${\bf Z}_X$ the constant sheaves of the integer group on $M$ and
$X$. Since $X$ is simply connected,  $$R^i f_* {\bf Z}_M \cong {\bf
Z}_X^{\oplus b_i(F)},$$ where $b_i(F)$ denotes the $i$-th Betti
number of a fiber $F$. From the beginning of the Leray spectral
sequence for $f$, we get an exact sequence
$$ H^1({\bf Z}_M) \longrightarrow H^0(R^1f_* {\bf Z}_M)
\longrightarrow H^2(f_* {\bf Z}_M).$$ The first term vanishes
because $M$ is simply connected. Consequently, $b_1(F) \leq b_2(X)$,
which is absurd because $b_1(F) = 2n$ and $b_2(X) =1$.
 $\Box$

\medskip
Now let $X$ be as in Theorem 1.2. By Theorem 1.1,
 the base manifold $X$ is a Fano manifold with $b_2(X)= b_4(X)=1$ and $f$ is a
 Lagrangian fibration.
 Fix a minimal component ${\cal K}$ of $X$ and let $C$ be a general
member of ${\cal K}.$ We want to apply Proposition 3.4 to the
complete curve $C \subset X$ and the Lagrangian fibration $f:M
\rightarrow X$. However, we have little information about $H^0(C,
T^*(X))$. What we have is the information on $H^0(\hat{C},
\nu^*T^*(X))$ where $\nu: \hat{C} \rightarrow C$ is the
normalization. To remedy this, we lift the Lagrangian fibration
$f$ to the normalization as follows.

\medskip
Recall that  $\nu: \hat{C} \rightarrow C \subset X$ is an
immersion of $\P_1$.
 By analytic
continuation, we can find an embedding of $\hat{C}$ in a complex
manifold $B_C$ of dimension $n$ and a holomorphic map $\alpha: B_C
\rightarrow M$ that is unramified at each point of $B_C$. The germ
of $B_C$ along $\hat{C}$ is uniquely determined by the germ of $X$
along $C$. By pulling back the Lagrangian fibration $f: M
\rightarrow X$ by the unramified holomorphic map $\alpha: B_C
\rightarrow X$, we get a Lagrangian fibration $$f_C:= \alpha^*f:
M_C \longrightarrow B_C.$$ Now we apply Proposition 3.4 to $f_C$
and the complete curve $\hat{C} \subset B_C$. We have an abelian
variety $A_{\hat{C}}$ acting on $f_C^{-1}(\hat{C})_{\rm red}$ in
the notation of Proposition 3.4.

\medskip
{\bf Proposition 4.2} {\it Let $f:M \rightarrow X$ be as in
Theorem 1.2 and ${\cal K}$ be a minimal component on $X$. Let
$\nu: \hat{C} \rightarrow C \subset X$ be the normalization of a
general member of ${\cal K}$ and $\alpha: B_C \rightarrow X$ be
the associated unramified holomorphic map as explained above. Let
$f_C : M_C \rightarrow B_C$ be the pull-back of $f$ by $\nu$. Then

 (1) $\dim A_{\hat C} \leq n-1$, and

 (2) for a point $b \in \hat{C}$ outside the critical set  and a
 point $z \in f_C^{-1}(b),$ the tangent space $T_z(A_{\hat{C}} \cdot z)$
 contains the $\iota_{b,z}$-image of the evaluation of $H^0(\hat{C}, T^*(B_C)) = H^0(\hat{C},
 \nu^*T^*(X))$ at the point $b$. }

 \medskip
 {\it Proof}. Note that the critical set $D$ of $f$  is an  ample hypersurface
  as a consequence of the condition $b_2(X)=1$ and Proposition 4.1.
   Thus $C$ intersects the critical set $D$
of
 $f$ and $\hat{C}$ intersects the critical set of $f_C$. Thus (1) and (2)
 follow from Proposition 3.4 (iii) and (iv).
 $\Box$

\medskip
{\bf Proposition 4.3} {\it Let $f: M \rightarrow X$ be as in
Theorem 1.2. For any minimal component ${\cal K}$ on $X$, the
variety of minimal rational tangents ${\cal C}_x$ for a general
point $x \in X$ is linear, i.e., each of its components is a
linear subspace in ${\P}T_x(X).$}

\medskip
{\it Proof}. Fix a general point $x \in X$ and a component ${\cal
C}_1$ of the variety of minimal rational tangents ${\cal C}_x
\subset \P T_x(X)$. Let ${\cal K}_1$ be the corresponding
component of ${\cal K}_x$. For a general member $C$ of ${\cal
K}_1$, $C$ is non-singular at $x$. Choose a point $z \in
f_C^{-1}(x).$ For simplicity, we will identify $f_C^{-1}(x)$ and
$f^{-1}(x)$ by the obvious isomorphism.  By Proposition 4.2 (2),
the image in $T^*_x(X)$ of the evaluation of $H^0(\hat{C},
\nu^*T^*(X))$ at $x$ is contained in $
\iota_{x,z}^{-1}(T_z(A_{\hat{C}} \cdot z))$ in the notation of
Proposition 3.4 (iv). Now let us vary  the general member $C$ of
${\cal K}_1$. The abelian subvariety $A_{\hat{C}}\cdot z $ of the
abelian variety $f^{-1}(x)$ with $z$ as the origin remains
unchanged, because there is no continuous family of abelian
subvarieties in a fixed abelian variety. Thus the image of the
evaluation of $H^0(\hat{C}, \nu^*T(X))$ is contained in a fixed
linear subspace $\iota_{x,z}^{-1} (T_z(A_{\hat{C}} \cdot z))$.
This linear subspace has dimension $<n$ by Proposition 4.2 (1).
 Thus by Proposition 2.1, the dual variety of the component
${\cal C}_1$ of ${\cal C}_x$ in ${\bf P}T^*_x(X)$ is degenerate.
Now apply Proposition 2.2  to complete the proof. $\Box$

\medskip
By Proposition 4.3, to prove Theorem 1.2,  we may assume that $X$
satisfies (Assumption) in Section 2. In particular, we have the
family of $\P_{p+1}$'s described in Proposition 2.4 and apply the
results in Section 2. The main consequence is the following.

\medskip
{\bf Proposition 4.4} {\it Let $H \subset X$ be an irreducible
hypersurface  such that a minimal component ${\cal K}$  is
multivalent on $H$. Then $H$ is not contained in the critical set
$D$ of the Lagrangian fibration $f: M \rightarrow X$. }

\medskip To prove Proposition 4.4, we need two  elementary
lemmata. First, some notation.
 Suppose that we are given an abelian variety $A$ acting on another
 abelian variety
$A'$ by a morphism $ \beta: A \times A' \rightarrow A'.$ Fixing a
point $z \in A'$, the orbit map $ \beta( \cdot, z): A \rightarrow
A'$ induces a homomorphism $\beta_*: H_1(A, {\bf Z}) \rightarrow
H_1(A', {\bf Z})$. In fact, this homomorphism does not depend on
the choice of $z$, because the translation on $A'$ acts trivially
on $H_1(A', {\bf Z})$.   We have the following two lemmata.

\medskip
{\bf Lemma 4.5} {\it Let $A_1$ and $A_2$ be two abelian varieties
acting on an abelian variety $A'$, by  two morphisms $\beta_1: A_1
\times A' \rightarrow A'$ and $\beta_2: A_2 \times A' \rightarrow
A'$. Suppose that for a point $z \in A'$, the tangents to the
orbits $T_z(A_1 \cdot z)$ and $T_z(A_2 \cdot z)$ span $T_z(A')$.
Then the images of the two homomorphisms
$$\beta_{1*}: H_1(A_1, {\bf Z}) \longrightarrow H_1(A', {\bf Z})
\mbox{ and } \beta_{2 *}: H_1(A_2, {\bf Z}) \longrightarrow
H_1(A', {\bf Z})$$ generate a subgroup of finite index in $H_1(A',
{\bf Z})$.}

\medskip
{\it Proof}. $A'_1:=A_1 \cdot z$ and $A'_2:=A_2 \cdot z$ are
abelian subvarieties in $A'$. Since $T_z(A'_1)$ and $T_z(A'_2)$
span $T_z(A')$,  there is a basis of $T_z(A')$ consisting of
lattice vectors of $A'_1$ and $A'_2$ under the universal covering
map $T_z(A') \rightarrow A'$. Since the images of $\beta_{1*}$
(resp. $\beta_{2*}$) is of finite index in $H_1(A'_1, {\bf Z})$
(resp. $H_1(A'_2, {\bf Z}))$,  Lemma 4.5 follows. $\Box$

\medskip
{\bf Lemma 4.6} {\it Let $f: {\cal A} \rightarrow B$ be a smooth
abelian fibration. Assume that there exists an effective fiberwise
action of an  abelian variety $A$ on ${\cal A}$ given by  $\beta:
A \times {\cal A} \rightarrow {\cal A}$. For any $b \in B$, let
$\beta (b): A \times f^{-1}(b) \rightarrow f^{-1}(b)$ be the
action on the fiber and let
$$\beta(b)_*: H_1(A, {\bf Z}) \longrightarrow H_1(f^{-1}(b), {\bf Z})$$
be the  homomorphism induced by $\beta(b)$. Then the image of
$\beta(b)_*$ is fixed under the monodromy action of $\pi_1(B,b) $
on $H_1(f^{-1}(b), {\bf Z})$. }

\medskip
{\it Proof}. As $b$ varies over $B$, $\beta(b)_*$ defines a
homomorphism from the constant  system $H_1(A, {\bf Z})$ on $B$
into the local system defined by $H_1(f^{-1}(b), {\bf Z})$ on $B$.
Thus the image is a constant system,  unchanged under the
monodromy. $\Box$

\medskip
{\it Proof of Proposition 4.4}.  Let us use the notation of
Proposition 2.6. Let $y \in H$ be a general point of $H$ and
choose neighborhoods $W \subset W_0$ as in Proposition 2.6. For a
general point $x \in W \setminus H$, we get $\zeta_1, \zeta_2 \in
Z_o$ satisfying the properties in Proposition 2.6.  For each $i
=1, 2$, choose a neighborhood $U_i$ of $\zeta_i$ in $Z$ such that
$\mu$ is unramified on $B_i:=\rho^{-1}(U_i)$. Pulling back $f: M
\rightarrow X$ by the unramified map $\mu|_{B_i}$, we have a
Lagrangian fibration $f_i: M_i \rightarrow B_i$. Applying
Proposition 3.4 to $f_i$ and the complete variety
$\rho^{-1}(\zeta_i)$, we get an action of an abelian variety $A_i$
on $f_i^{-1}(\rho^{-1}(\zeta_i))_{\rm red}$. Note that the loop
$\gamma_i$ on the smooth locus of $P_{\zeta_i}$  is naturally
lifted to a loop in $\rho^{-1}(\zeta_i)$, which we denote by the
same symbol $\gamma_i$. Pick a point $z \in f^{-1}(x)$. Let $x_i
\in \rho^{-1}(\zeta_i)$ be the point over $x$ and $z_i \in
f_i^{-1}(x_i)$ be a point on the fiber corresponding to $z$.  The
monodromy of the smooth Lagrangian fibration along $\gamma_i$
fixes the image of $H_1(A_i, {\bf Z})$ in $$H_1(f^{-1}_i(x_i),
{\bf Z}) = H_1(f^{-1}(x), {\bf Z})$$ by Lemma 4.6.  But the two
subspaces $\iota_{x,z}^{-1}(T_z(A_1 \cdot z))$ and
$\iota_{x,z}^{-1}(T_z(A_2 \cdot z))$ span the whole $T^*_x(X)$ by
Proposition 2.5 and Proposition 4.2 (2). It follows that the
monodromy action of $\pi_1(W_0 \setminus H)$ on $H_1(f^{-1}(x),
{\bf Z})$ is trivial by Lemma 4.5. Thus the hypersurface is not in
the critical set by Proposition 3.2. $\Box$

\medskip
An immediate consequence of Proposition 2.8, Proposition 4.1 and
Proposition 4.4 is the following.

\medskip
{\bf Proposition 4.7} {\it  Let $f: M \rightarrow X$ be as in
Theorem 1.2. Suppose $X$ is different from $\P_n$. Then for any
minimal component ${\cal K}$, $p= \dim {\cal C}_x =0$.}

\section{Use of the integrability of the distribution defined by a pair of rational curves}

By Proposition 4.7, to prove Theorem 1.2,  we may assume that the
base manifold $X$ has a minimal component ${\cal K}$ with $p=0$.
We will make this assumption throughout this section. We  start
with the following observation.

\medskip
{\bf Proposition 5.1} {\it Let $X$ be as above with $p=0$. Let
$\mu: X' \rightarrow X$ be as in Proposition 2.4. If $E \subset X$
is a component of the branch locus of $\mu$, then $E$ is contained
in the critical set $D$ of $f:M \rightarrow X$.}

\medskip
{\it Proof}. Suppose not. Let $y$ be a general point of $E$ outside
the critical set $D$. Let $R$ be a component of the ramification
locus of $\mu$ such that $\mu(R) = E$. Let $z \in R$ be a point in
$\mu^{-1}(y)$. We can choose an open neighborhood $W'$ of $z$ and an
open neighborhood $W$ of $y$ with $W \cap D = \emptyset$ such that

(i) $W$ and $W'$ are biholomorphic to a polydisc with $E \cap W
\subset W$ and $R \cap W' \subset W'$ biholomorphic to the
coordinate hyperplane;

(ii) there exists a Lagrangian section $\Sigma \subset f^{-1}(W)$
of $f$ over $W$;

 (iii) $W' \setminus R \subset \rho^{-1}(Z_o)$; and

(iv) $\mu|_{W'}: W' \rightarrow W$ is a cyclic branched covering of
degree $>1$.

Pulling back $f$ by the unramified holomorphic map $\mu|_{W'
\setminus R}$, we get a smooth Lagrangian fibration $f': M'
\rightarrow W' \setminus R$ with a Lagrangian section $\Sigma'
\subset M'$. Applying Proposition 3.7, we get a family of abelian
subvarieties of dimension $n-1$, say ${\cal A} \subset M'$, over $W'
\setminus R$. By Proposition 2.4 (c), we know that the fibers of
$\rho$ on $W' \setminus R$ are sent to $W \setminus E$ to define a
multi-valued foliation. The relation between the family of  abelian
subvarieties ${\cal A}$ and the fibers of $\rho$ in Proposition 3.7
implies that the image of ${\cal A}$ in $f^{-1}(W \setminus E)$ is a
multi-valued family of abelian subvarieties. Thus the Galois group
of the cyclic covering $\mu|_{W'}$ acts non-trivially on the family
of lattices for ${\cal A}$. This means that the monodromy of
$\pi_1(W \setminus E)$ acts non-trivially on $H_1(f^{-1}(b), {\bf
Z})$ for a point $b \in W \setminus E$. This is absurd because the
Lagrangian fibration $f$ is smooth over $W$. $\Box$

\medskip
Now we introduce a new notion.
 The family
$Z$ in Proposition 2.4 is a ${\bf P}_1$-bundle over an open subset
$Z_o$ of $Z$. Let $x \in X$ be a general point so that all members
of ${\cal K}_x$ are smooth at $x$ and have distinct tangent
vectors at $x$ by Proposition 2.5. Furthermore, we can assume that
$x$ lies
 outside $\mu(\rho^{-1}(Z\setminus Z_o))$ and the branch locus of
$\mu$. Let $\zeta_1, \zeta_2 \in Z_o$ be two distinct points with
$x \in P_{\zeta_1} \cap P_{\zeta_2}$. Let $y_1 =
\rho^{-1}(\zeta_1) \cap \mu^{-1}(x)$ and $y_2 = \rho^{-1}(\zeta_2)
\cap \mu^{-1}(x)$. Choose a small neighborhood $W$ of $x$, $W_1$
of $y_1$ and $W_2$ of $y_2$ such that

(i)  $W_1$, $W_2$ and $W$ are biholomorphic to one another by
$\mu$;

(ii) fibers of $\rho$ in $W_1$ and fibers of $\rho$ in $W_2$
induce two transversal holomorphic foliations of rank 1 on $W$.

The span of these two foliations define a distribution of rank 2
on $W$, which we denote by ${\cal D}$. We call ${\cal D}$ the {\it
distribution defined by the pair $P_{\zeta_1}$ and $P_{\zeta_2}$}.

\medskip
{\bf Proposition 5.2} {\it Let $X$ be a Fano manifold with
$b_2(X)=1$ having a minimal component ${\cal K}$ with $p=0$. Let $x
\in X$ be a general point and $C_1, C_2$ be two distinct members of
${\cal K}_x$. If $X$ is a base manifold of a fibration $f:M
\rightarrow X$ as in  Theorem 1.2, then the distribution ${\cal D}$
defined by the pair $C_1$ and $C_2$ in a neighborhood of $x$ is
integrable.}

\medskip
{\it Proof}. Let $B := \rho^{-1}(Z_o)$, then $\mu$ is unramified
on $B$. Pulling back $f:M \rightarrow X$ by $\mu$, we get a
Lagrangian fibration $g: M' \rightarrow B$ and a smooth fibration
by curves $\rho: B \rightarrow Z_o$.  Let $W$ be a neighborhood of
$x$ and $W_1, W_2 \subset B$ be the open subsets used in the
definition  of ${\cal D}$. By shrinking $W$ if necessary, we may
assume that $f$ has a Lagrangian section over $W$, hence the
induced Lagrangian fibrations over $W_1$ and $W_2$ also have
Lagrangian sections. Applying Proposition 3.7 to $W_1$ and $W_2$,
we get  two families ${\cal A}_1$ and ${\cal A}_2$ of abelian
subvarieties of dimension $n-1$ inside the Lagrangian fibration
$f|_{f^{-1}(W)}$. The annihilator ${\cal D}^{\perp} \subset
T^*(W)$ of the distribution ${\cal D} \subset T(W)$ corresponds,
in the sense of Proposition 3.6,  to the family of abelian
subvarieties of dimension $n-2$ defined by the intersection ${\cal
A}_1 \cap {\cal A}_2$. Thus ${\cal D}$ is integrable by
Proposition 3.6. $\Box$

\medskip
{\bf Proposition 5.3} {\it Let us assume the situation of
Proposition 5.2 and use the notation of Proposition 2.4.  Given a
general member $C$ of ${\cal K}$, let $C'$ be a component of
$\mu^{-1}(C)$ such that
$$\mu|_{C'}: C' \longrightarrow C$$ is finite of degree $>1$, which exists by Proposition
2.3. Let  $$h: \hat{C}' \longrightarrow \widehat{\rho(C')}$$ be
the lift of $$\rho|_{C'}: C' \longrightarrow \rho(C')$$ to the
normalizations of $C'$ and $\rho(C')$. Then $h$ has a ramification
point $z \in \hat{C}'$ such that the image of $h(z)$ in $\rho(C')$
lies in $Z_o$.}

\medskip
We  need two lemmata.

\medskip
{\bf Lemma 5.4} {\it  Let $E$ be a component of the branch divisor
of $\mu$. Let $E'$ be an irreducible component of $\mu^{-1}(E)$
which is dominant over ${\cal K}$ by $\rho$. Then $C'$ is disjoint
from $E'$.}

\medskip {\it Proof}. Otherwise, through a general point $x'$ of $E'$,
  we have two distinct curves, $C'$ and a fiber of
$\rho$, neither of which are contained in $E'$. Since $\mu$ is
unramified at $x'$ by Proposition 2.4, the images of these curves
under $\mu$ are of the form $P_{\zeta_1}, P_{\zeta_2}$ with
$\zeta_1 \neq \zeta_2$. Since these two curves pass through $x =
\mu(x')$, which  is a general point of $E$, ${\cal K}$ is
multivalent on $E$. Then $E$ is not in the critical set of $f$ by
Proposition 4.4,  a contradiction to Proposition 5.1. $\Box$

\medskip
{\bf Lemma 5.5} {\it There exist a family of members of ${\cal K}$
$$\{C_t, t \in \Delta, C=C_0 \} $$ and the associated deformation
$$\{ C'_t, t \in \Delta, C'=C'_0\}$$  such that for each $t \in
\Delta$,

(i) $C'_t$ is a component of $\mu^{-1}(C_t)$;

(ii)  $\mu|_{C'_t}: C'_t \rightarrow C_t$ is finite of degree $>1$
over $C_t$;

(iii) $\rho(C'_t) = \rho(C')$.}

\medskip
{\it Proof}.  Let $x \in C$ be a general point. Pick a point $x' \in
C'$ with $\mu(x') = x$ and set $C^{x} := \mu(\rho^{-1}(\rho(x')))$.
Then $C^x$ is a member of ${\cal K}_x$. The triviality of the normal
bundle of $C$ implies that there exists a `deformation of $C$ along
$C^x$', namely, a unique (up to reparametrization) deformation $\{
C_t, t \in \Delta\}$ of $C=C_0$ satisfying $C_t \cap C^x \neq
\emptyset$.
 More precisely, there
exists a unique component $C^{\sharp}$ of $\mu^{-1}(C^x)$
different from $\rho^{-1}(\rho(x'))$ that contains the point
$\rho^{-1}(\zeta) \cap \mu^{-1}(x)$ where $P_{\zeta} =C$. Then the
germ of $\rho(C^{\sharp})$ near $\zeta$ gives a deformation $\{
C_t, t \in \Delta, C=C_0 \}$ such that each $C_t$ belongs to
${\cal K}$ and $C_t \cap C^{x} \neq \emptyset$. The germ of such a
deformation of $C$ is uniquely determined by $C^x$, up to
reparametrization.

By the generality of the choice of $C$, we have the associated
deformation $\{ C'_t, t \in \Delta\}$ of $C'= C'_0$ such that
$\mu|_{C'_t}: C'_t \rightarrow C_t$ is finite of degree $>1$ for any
$t \in \Delta$.

It remains to show that $\rho(C'_t) = \rho(C')$ for each $t \in
\Delta$. This follows from the fact that the leaf through $x$ of the
distribution ${\cal D}$ defined by $C$ and $C^x$  gives  the germ of
the surface traced out by  the deformation of $C$ along $C^x$ and at
the same time the germ of the surface traced out by the deformation
of $C^x$ along $C$.  More precisely, it suffices to show that
$$ (*) \;\; \rho^{-1}(z) \cap C'_t \neq \emptyset \mbox{ for each
} z \in \rho(C') \mbox{ and each } t \in \Delta.$$  Clearly
$\rho^{-1}(\rho(x')) \cap C'_t \neq \emptyset.$ Suppose that we
vary the choice of $x \in C$ in the above construction of $C_t$,
say $\{x_s \in C, s \in \Delta, x_0 =x\}$. Let $x'_s \in C'$ be
the associated variation of $x'=x_0$ and $C^{x_s}:= \mu(
\rho^{-1}(\rho(x_s))$. By Proposition 5.2, the deformation of $C$
along $C^{x_s}$ are, up to reparametrization, the same as the
deformation of $C$ along $C^x$ for any $s \in \Delta$. Thus
$\rho^{-1}(\rho(x'_s)) \cap C'_t \neq \emptyset$ for each $s \in
\Delta$. This proves $(*)$ for $z$ in a neighborhood of
$\rho(x')$. It follows that $(*)$ holds for all $z \in \rho(C')$.
$\Box$

\medskip
{\it Proof of Proposition 5.3}.
 Suppose  that $h$ is unramified over $\rho(C') \cap Z_o$.
  Let us use the deformation $C_t$ constructed in
 Lemma 5.5. By the
generality of $C$, we may assume that for each $t \in \Delta$ the
holomorphic map $$h_t: \hat{C}'_t \rightarrow \widehat{\rho(C'_t)} =
\widehat{\rho(C')}, \;\;  h_0 =h,$$ which is the lift of
$\rho|_{C'_t}$ to the  normalization of curves, is unramified over
$\rho(C') \cap Z_o$. Since $h_t$ is a continuous family of coverings
of the Riemann surface $\widehat{\rho(C')}$ with fixed branch locus,
we can find a biholomorphic map
$$ (\clubsuit) \;\;\; \psi_t: \hat{C}' \rightarrow \hat{C}'_t, \;\; \psi_0 = {\rm Id}_{\hat{C}'} \mbox{ with } h = h_t \circ
\psi_t,$$ which depends holomorphically on $t$ (e.g. [Sh, p. 32,
Corollary 1].).

From Lemma 5.5 (ii), there are at least two distinct points in
$\hat{C}_t$, say $a_t \neq b_t \in \hat{C}_t$, such that the
corresponding points in $C_t$ lie in the branch divisor of $\mu$ in
$X$. Let $\{ 0, \infty\} \subset \P_1$ be two distinct points on the
projective line. We can choose a family of biholomorphic maps $\{
\sigma_t: \hat{C}_t \rightarrow \P_1, \; t \in \Delta \}$ such that
$\sigma_t(a_t) = 0$ and $ \sigma_t(b_t) = \infty$ for each $t \in
\Delta$.  Denote by
 $\mu_t: \hat{C}'_t \rightarrow \hat{C}_t$  the lift of
$\mu|_{C'_t}$ to the normalization of curves. Then
 $$ \{ \varphi_t: \hat{C}' \longrightarrow
\P_1, \;\;\; \varphi_t := \sigma_t \circ \mu_t \circ \psi_t, \; t
\in \Delta \}$$ is a family of meromorphic functions on the compact
Riemann surface $\hat{C}'$.

By Lemma 5.4, for each component $E$ of the branch divisor of $\mu$,
the intersection of $C'_t$ with $\mu^{-1}(E)$ has a fixed image in
$\rho(C')=\rho(C'_t)$, independent of $t \in \Delta$. This implies
that there is a finite subset $Q \subset \widehat{\rho(C')}$,
independent of $t$, such that
$$\mu_t^{-1}(a_t) \cup \mu_t^{-1}(b_t) \subset h_t^{-1}(Q)$$ for any $t
\in \Delta$. Then $$ \varphi_t^{-1}(0) = \psi^{-1}_t \circ
\mu_t^{-1} \circ \sigma_t^{-1}(0) = \psi^{-1}_t(\mu_t^{-1}(a_t))
\subset  \psi_t^{-1}(h_t^{-1}(Q)) $$ for all $t \in \Delta$. Since
$\psi_t^{-1}(h_t^{-1}(Q))= h^{-1}(Q)$ by the choice of $\psi_t$ in
$(\clubsuit)$, $ \varphi_t^{-1}(0) \subset h^{-1}(Q)$ for any $t \in
\Delta$. Consequently,  $\varphi_t^{-1}(0) = \varphi_0^{-1}(0)$ for
all $t \in \Delta$. By the same argument we get
$\varphi_t^{-1}(\infty) = \varphi_0^{-1}(\infty)$ for all $t \in
\Delta$. In other words, the family of meromorphic functions
$\varphi_t$ have the same zeroes and the same poles on the Riemann
surface $\hat{C}'$. This implies that for any $z \in \P_1$ and $t
\in \Delta$, $\varphi^{-1}_t(z) = \varphi^{-1}_0(z)$. It follows
that for any $w_1, w_2 \in \hat{C}'$ and any $t \in \Delta$, $$
(\diamondsuit) \;\;\;\; \varphi_t(w_1) = \varphi_t(w_2) \mbox{ if
and only if } \varphi_0(w_1) = \varphi_0(w_2).$$

Since $\mu|_{C'}$ is  finite of degree $>1$ by our assumption, we
can choose two points $\alpha \neq \beta \in \hat{C}'$ such that
$\varphi_0(\alpha) = \varphi_0(\beta).$  Furthermore, denoting by
$\bar{\alpha}\in \rho(C')$ (resp. $\bar{\beta} \in \rho(C')$)  the
point corresponding to $h_0(\alpha) \in \widehat{\rho(C')}$ (resp.
$h_0(\beta)\in \widehat{\rho(C')}$) under the normalization, we may
assume that $$ (\heartsuit) \;\;\; \bar{\alpha} \mbox{ and }
\bar{\beta} \mbox{ are two distinct points in } Z_o.$$ From
$(\diamondsuit)$, we have  $\varphi_t(\alpha) = \varphi_t(\beta)$
for all $t \in \Delta$. Since $\varphi_t = \sigma_t \circ \mu_t
\circ \psi_t$ and $\sigma_t$ is biholomorphic, we see that $$
(\spadesuit) \;\;\;\;  \mu_t \circ \psi_t(\alpha) = \mu_t \circ
\psi_t(\beta) \mbox{ for all } t \in \Delta.$$

 Denote by $$ \alpha_t \in C'_t \subset X' \;\;\; \mbox{ (resp. }
 \beta_t \in C'_t \subset X') $$ the point corresponding to
  $\psi_t(\alpha) \in \hat{C}'_t$ (resp.
 $\psi_t(\beta) \in \hat{C}'_t$)  under the normalization. Then the
 locus $$A := \{ \alpha_t \in X', t \in \Delta \} \;\;\; \mbox{ (resp.
  } B := \{ \beta_t \in X', t \in \Delta \})$$
covers a non-empty  open subset in the fibre
$\rho^{-1}(\bar{\alpha})$ (resp. $\rho^{-1}(\bar{\beta})$). Thus
$\mu(A)$ (resp. $\mu(B)$) covers a non-empty open subset in $$
P_{\bar{\alpha}} := \mu(\rho^{-1}(\bar{\alpha})) \;\;\; \mbox{
(resp. } P_{\bar{\beta}} :=\mu( \rho^{-1}(\bar{\beta} ))).$$ Since
$\mu(A)$ (resp. $\mu(B)$) is the locus of points corresponding to
$\mu_t \circ \psi_t (\alpha)$ (resp. $\mu_t \circ \psi_t (\beta)$)
by the normalization $\hat{C}_t \rightarrow C_t$, the equality
$(\spadesuit)$ above implies that $\mu(A) = \mu(B)$. Consequently,
$$P_{\bar{\alpha}}= P_{\bar{\beta}},$$ a contradiction to
$(\heartsuit)$ and  Proposition 2.4 (c). $\Box$

\medskip
{\bf Proposition 5.6} {\it In the situation of Proposition 5.3,
there exists an irreducible hypersurface $H \subset X$ such that
through a general point $x \in H$, there are two distinct members
$C_1$ and $C_2$ of ${\cal K}$ with $T_x(C_1) = T_x(C_2)$ in
$T_x(X)$ where $T_x(C_1)$ (resp. $T_x(C_2)$) denotes the tangent
space at $x$ of a component of the germ of $C_1$ (resp. $C_2$) at
$x$. In particular, ${\cal K}$ is multivalent on $H$.}

\medskip
{\it Proof}. In the situation of Proposition 5.3, let $z \in C'$ be
the image of a ramification point of $h$ such that $\rho(z) \in
Z_o$. Then  $\mu$ is unramified at $z$ and the curve
$\mu(\rho^{-1}(\rho(z)))$ in $X$ is an immersed $\P_1$. It follows
that $C'$ is immersed at $z$ and one of the component of the germ of
$C'$ at $z$ must be tangent to $\rho^{-1}(\rho(z))$ because $z$ is a
ramification point of $h$. But then
  the two members  $C$ and $\mu(\rho^{-1}(\rho(z)))$ of ${\cal K}$
  have non-empty
intersection at $\mu(z)$, sharing a common tangent. As $C$ varies,
this intersection point also varies to define a hypersurface $H$.
This proves Proposition 5.6. $\Box$

\medskip
Finally, the next proposition leads to a contradiction with
Proposition 4.4 and  Proposition 5.6, completing the proof of
Theorem 1.2.

\medskip
{\bf Proposition 5.7} {\it The hypersurface $H$ in Proposition 5.6
is contained in the critical set $D$ of $f$.}

\medskip
{\it Proof}. Suppose not.  Let $y\in H$ be a general point. By the
definition of $H$, there are two points $y_1, y_2
 \in \mu^{-1}(y)$ where $\mu$ is unramified,
 such that $P_{\zeta_1}, \zeta_1 = \rho(y_1),$ and $P_{\zeta_2}, \zeta_2=
 \rho(y_2)$ are two distinct members of ${\cal K}$ tangent to each other at $y$.
Choose a neighborhood $W$ of $y$ and a neighborhood $W_i$ of $y_i$
for each $ i= 1,2$ such that $W$, $W_1$ and $W_2$ are all
biholomorphic by $\mu$. By shrinking $W$ if necessary, we can assume
that $f|_{f^{-1}(W)}$ is a smooth Lagrangian fibration with a
Lagrangian section over $W$. The same holds for the pull-backs of
$f$ to $W_1$ and $W_2$.  Thus by Proposition 3.7, we get two
families of abelian subvarieties of dimension $n-1$ in the
Lagrangian fibration $f^{-1}(W) \rightarrow W$ such that each
abelian subvariety in the fiber is tangent to the image of the
conormal bundle of deformations of $C$ and $C'$. Since $C$ and $C'$
are tangent along $H$, the two families of abelian subvarieties
coincide along $H$, but not at a general point of $W$. This is a
contradiction to the following elementary lemma. $\Box$

\medskip
{\bf Lemma 5.8} {\it Given a family of abelian variety ${\cal
A}\rightarrow S$ and two families of abelian subvarieties ${\cal
A}_1 \rightarrow S$ and ${\cal A}_2 \rightarrow S$, if the fibers
of ${\cal A}_1$ and ${\cal A}_2$ coincide at some point of $S$,
then ${\cal A}_1 = {\cal A}_2$. }

\medskip
{\it Proof}. Consider the family of quotient abelian varieties
${\cal A}/{\cal A}_1 \rightarrow S$ and the relative group quotient
projection ${\cal A} \rightarrow {\cal A}/{\cal A}_1$. This
projection sends one fiber of ${\cal A}_2$ to a point. So it must
send each fiber of ${\cal A}_2$ to a point by [Mu, Proposition 6.1].
Thus ${\cal A}_1 = {\cal A}_2$. $\Box$

\medskip
{\bf Acknowledgment} It is a pleasure to thank  Yasunari Nagai,
Yong-Geun Oh, and Justin Sawon  for beneficial discussions. I am
particularly grateful to Daisuke Matsushita, who kindly sent me his
preprint [Ma4] and  explained several points of his works to me. A
very special thank goes to Keiji Oguiso, for many discussions on
holomorphic symplectic geometry and very valuable comments on the
first draft of this paper.

\bigskip
\bigskip
{\bf References}

\medskip

[CMS] Cho, K., Miyaoka, Y. and Shepherd-Barron, N.:
Characterizations of projective space and applications to complex
symplectic manifolds. In {\it Higher dimensional birational
geometry (Kyoto, 1997)}, Adv. Stud. Pure Math. {\bf 35} (2002)
1--88

 [GuSt] Guillemin, V. and Sternberg, S.: {\it Symplectic
techniques in physics}. Second edition. Cambridge University
Press, Cambridge, 1990

[Hu] Huybrechts, D.: Compact hyperk\"ahler manifolds. In {\it
Calabi-Yau manifolds and related geometries (Nordfjordeid, 2001)},
161--225, Universitext, Springer, Berlin, 2003

 [Hw] Hwang, J.-M.: Deformation of
holomorphic maps onto Fano manifolds of second and fourth Betti
numbers 1. Ann. Inst. Fourier {\bf 57} (2007) 815-823

[HwMo1] Hwang, J.-M. and Mok, N.:  Cartan-Fubini type extension of
holomorphic maps for Fano manifolds of Picard number 1. Journal
Math. Pures Appl. {\bf 80}  (2001) 563-575

 [HwMo2] Hwang, J.-M. and Mok,
N.: Birationality of the tangent map for minimal rational curves.
Asian J. Math. {\bf 8}, {\it Special issue dedicated to Yum-Tong
Siu}, (2004) 51-64

[HwOg] Hwang, J.-M. and Oguiso, K.: Characteristic foliation on
the discriminantal hypersurface of a holomorphic Lagrangian
fibration. preprint (arXiv:0710.2376)

 [HwRa] Hwang, J.-M. and Ramanan,
S.: Hecke curves and Hitchin discriminant. Ann. scient. Ec. Norm.
Sup. {\bf 37} (2004) 801-817

 [Ka] Kaup, W.:
Infinitesimale Transformationsgruppen komplexer R\"aume. Math.
Annalen {\bf 160} (1965) 72-92

 [Ma1] Matsushita, D.: On fibre space structures of a projective
  irreducible symplectic manifold. Topology {\bf 38} (1999),
  79--83,
  and addendum: Topology {\bf 40} (2001) 431--432

[Ma2] Matsushita, D.: Equidimensionality of Lagrangian fibrations on
holomorphic symplectic manifolds. Math. Res. Lett. {\bf 7} (2000)
389-391

 [Ma3] Matsushita, D.: Higher direct images of
dualizing sheaves of Lagrangian fibrations. Amer. J. Math. {\bf 127}
(2005) 243--259

 [Ma4] Matsushita, D.:  A canonical bundle formula of projective
Lagrangian fibrations. preprint (arXiv:0710.0122)

 [Mu] Mumford, D. and Fogarty, J.: {\it Geometric invariant
theory}. Second enlarged edition. Springer Verlag, 1982

[Ng] Nagai, Y.: Dual fibration of a projective Lagrangian
fibration. preprint

 [Nk] Nakayama, N.: Compact K\"ahler manifolds whose universal
covering spaces are biholomorphic to ${\bf C}^n$. RIMS preprint
1230.

[Sh] Shokurov, V. V.: Riemann surfaces and algebraic curves. In
{\it Algebraic curves, algebraic manifolds and schemes}, Springer
Verlag, 1998

 \vspace{10mm}

Jun-Muk Hwang

Korea Institute for Advanced Study

207-43 Cheongnyangni-dong

Seoul 130-722, Korea

 jmhwang@kias.re.kr

\end{document}